\documentclass{amsart}
\usepackage{csquotes}
\usepackage{graphicx,amsaddr}
\usepackage{amsfonts}
\usepackage{amsmath}
\usepackage{mathtools}
\usepackage{amsthm}
\usepackage{amsfonts}
\usepackage{amssymb}
\usepackage{calc}
\usepackage{rotating}
\usepackage[all,cmtip]{xy}
\usepackage{pifont}
\usepackage{enumitem}
\usepackage{lscape}
\usepackage{array}

\theoremstyle{plain}
\newtheorem{defn}{Definition}[section]
\newtheorem{thm}[defn]{Theorem}

\newtheorem{prop}[defn]{Proposition}

\newtheorem{exam}[defn]{Example}
\newtheorem{conjec}[defn]{Conjecture}

\numberwithin{equation}{section}

\title[Spectral extension property]{The Spectral extension property in the unitization of Banach Algebras}

\author{H. V. Dedania, A. B. Patel*}

\address{Dept. of Mathematics, Sardar Patel University, Vallabh Vidyanagar 388120, Gujarat, India}
\email{hvdedania@gmail.com, avadhpatel663@gmail.com*}

\newcolumntype{C}[1]{>{\centering\arraybackslash}m{#1}}

\begin{document}



\begin{abstract}
Let $A$ be a non-unital Banach algebra and let $A_e = A \oplus {\mathbb C}1$ be the unitization of $A$. It is true that if $A_e$ has the spectral extension property (SEP), then $A$ has the same. Does the converse hold? In this paper, we give some necessary as well as some equivalent conditions.
\end{abstract}
\maketitle

\section{Introduction}

Let $A$ be a Banach algebra and let $a \in A$. Then $\sigma_A(a)$ and $r_A(a)$ denote the spectrum and the spectral radius of $a$ in $A$, respectively. These two concepts are extremely important in Banach Algebra Theory \cite[Section 5]{BoDu:73}. Let $A$ and $B$ be Banach algebras. Then $B$ is a \emph{supper Banach algebra} of $A$ if $A$ is a (not necessarily closed) subalgebra of $B$; the $B$ is a \emph{spectral extension} of $A$ if, further, $r_A(a) = r_B(a) \, (a \in A)$. The Banach algebra $A$ has the \emph{spectral extension property} (SEP) if every supper Banach algebra of $A$ is a spectral extension of $A$ \cite{Me:91}. Several equivalent conditions are given by several authors \cite{ToYo:89, BhDe:96}. A Banach algebra having SEP has a rich structure. Many strong results can be proved in such algebras \cite{Me:92}.

Let $A$ be a non-unital Banach algebra and let $A_e$ be its unitization. In genereal, any Banach algebra property of $A$ is preserved in $A_e$ also. For example: (i) $A$ is semisimple (resp. commutative, amenable, weakly amenable) iff $A_e$ is semisimple (resp. commutative, amenable, weakly amenable) \cite[Section 2.8]{Da:00}; (ii) $A$ has UUNP iff $A_e$ has UUNP \cite{DaDe:09}; (iii) In the case of Banach *-algebra, $A$ has UC*NP iff $A_e$ has UC*NP \cite{DeKa:13}.

Since $A$ is an ideal in $A_e$, it is easy to show that if $A_e$ has SEP, then $A$ has SEP. It is natural to ask a question, whether or not the converse holds true? Bhatt and Dedania answered this question in affirmative if $A$ is semisimple and commutative. Under some mild condition, the same was proved for semisimple Banach algebra in \cite{DeKa:13}.

In this paper, we prove that if $A$ has SEP, then both $A$ and $A_e$ have P-property and every norm on $A_e$ majorizes  a minimal norm. Moreover, we give some equivalent conditions. Surprisingly, one of them is the spatial numerical range (SNR) $V_A(a)$ of $a$ in $A$. The SNR is studied by several authors \cite{GaHu:89, Ta:00}. The SNR highly depends on the algebra $A$ and the norm on $A$. It is extensively studied in \cite{DePa:21}. Another equivalent condition is that, for any minimal norm $|\cdot|$ on $A_e$, the operator norm $|\cdot|_{op}$ and the $\ell^1$-norm $|\cdot|_1$ on $A_e$ are equivalent. In general, $|\cdot|_{op} \leq |\cdot|_1$ on $A_e$ but may not be equal \cite{DPa:21}.

\section{Main Results}
Let $(A, \|\cdot\|)$ be a normed algebra. Then $V_A(a)$ denotes the spatial numerical range (SNR) of $a$ in $(A, \|\cdot\|)$. We must note that the SNR depends on both $A$ and $\|\cdot\|$. It will be different even for two equivalent norms. It is studied in detail in \cite{DePa:21} with several examples. If $K \subset \mathbb C$, then $co(K)$ denotes the convex hull of $K$ and $\overline{co}(K)$ is the closure of $co(K)$. Let $\|\cdot\|$ be a norm on $A$. Then the operator norm $\|\cdot\|_{op}$ and the $\ell^1$-norm $\|\cdot\|_1$ on $A_e$ are defined as
\begin{eqnarray*}
\|a+\lambda 1\|_{op} & = & \sup\{\|a x + \lambda x \|: x \in A, \|x\| \leq 1\} \quad{(a + \lambda1 \in A_e)};\\
\|a+\lambda 1\|_1 & = & \|a\| + |\lambda| \quad{(a + \lambda1 \in A_e)}.
\end{eqnarray*}
It is clear that $\|\cdot\|_{op} \leq \|\cdot\|_1$ on $A_e$. The reader should refer \cite{DePa:21} for further properties.

The following result is a slight extension of the main theorem in~\cite{GaKo:93}; its proof is exactly the same. It will be used in the proof of our main theorem.

\begin{thm}\label{C}
Let $(A, \|\cdot\|)$ be a non-unital, normed algebra such that $\|\cdot\|$ is a regular norm on $A$ and $0 \in \overline{co}V_A(a) \; (a \in A)$. Then $$\|a + \lambda 1\|_{op} \leq \|a + \lambda 1\|_1 \leq 6 exp(1) \|a + \lambda 1\|_{op} \quad(\forall \lambda \in \mathbb C; \forall a \in A).$$
\end{thm}

Next result shows that the hypothesis of Theorem \ref{C} is satisfied under several mild conditions. A norm $|\cdot|$ on $A$ is called a \emph{Q-norm} if $r_A(a) \leq |a| \, (a \in A).$

\begin{prop}
Let $(A, \|\cdot\|)$ be a non-unital, normed algebra and $\|\cdot\|$ be a regular norm. If one of the following conditions holds, then $0 \in \overline{co}V_A(a) \, (a \in A)$; and hence $\|\cdot\|_{op} \cong \|\cdot\|_1$ on $A_e$ due to Theorem \ref{C}.
\begin{enumerate}
\item The norm $\|\cdot\|_{op}$ is a Q-norm on $A_e$;
\item $A$ is closed in $(A_e, \|\cdot\|_{op})$;
\item $\|\cdot\|$ is complete on $A$.
\end{enumerate}
\end{prop}

\begin{proof}
\emph{(1)} Define $\varphi_{\infty}(a + \lambda 1) = \lambda \; (a+ \lambda 1 \in A_e)$. We know that every complex homomorphism is continuous in any Q-norm. We have $\|\varphi_{\infty}\| = 1$ and $\varphi_{\infty}(1) = 1$. Therefore $0 = \varphi_{\infty}(a) \in V_{A_e}(a) \; (a \in A)$ with the norm $\|\cdot\|_{op}$. Therefore, by the \cite[Theorem 2.4(1)]{DePa:21}, we get $0 \in \overline{co}V_A(a) \; (a \in A)$.\\
\emph{(2)} $A$ is closed in $(A_e, \|\cdot\|_{op})$. Then the map $\varphi_{\infty}$ defined in the proof of statement \emph{(1)} above is continuous. Therefore, by the same arguments as in the proof of statement \emph{(1)} above, the result follows.\\
\emph{(3)} The norm $\|\cdot\|$ is complete and regular on $A$. Therefore $A = \ker \varphi$ is closed in $(A_e, \|\cdot\|_{op})$. Hence again the map $\varphi_{\infty}$ is continuous and so $0 \in V_A(a)$.
\end{proof}

\begin{exam}
(i) This example exhibits that Theorem \ref{C} is a generalization of \cite{GaKo:93}. Let $c_{00} = \{f: \mathbb N \longrightarrow \mathbb C: f(n) =0 \text{ eventually}\}$. Then $A = (c_{00}, \|\cdot\|_{\infty})$ is a non-complete, normed algebra. The Banach space dual $A^*$ of $A$ is linearly isomorphic to $(\ell^1, \|\cdot\|_1)$. Let $f = \sum_{n=1}^kf(n)\delta_n \in A$, let $g = \delta_{k+1} \in A$, and let $\varphi =\delta_{k+1} \in A^*$. Then $\|g\|_{\infty} = 1, \, \|\varphi\|_1 = \varphi(g) = 1$, and $0 = \varphi(fg) \in V_A(f)$.

(ii) It is not always true that $0 \in \overline{co} V_A(a)$. For example, if $A$ is any unital normed algebra with the identity $1_A$, then the spatial numerical range $V_A(1_A)$ of the identity $1_A$ is $\{1\}$.
\end{exam}

The results stated in the next proposition will be required in proofs. An algebra $A$ is \emph{faithful} if $a \in A$ and $aA = \{0\}$, then $a =0$. For example, every semisimple algebra is faithful (see \cite{BhDe:95}).

\begin{prop}\label{main1} Let $A$ be a non-unital, faithful algebra.
\begin{enumerate}
\item If $|\cdot|$ is a minimal norm on $A$, then $|\cdot|_{op} = |\cdot|$ on $A$;
\item If $|\cdot|$ is a minimal norm on $A_e$, then $|\cdot|_{op} = |\cdot|$ on $A_e$;
\item If $|\cdot|$ is a minimal norm on $A$, then $|\cdot|_{op}$ is a minimal norm on $A_e$;
\item If $A$ is dense in $(A_e, |\cdot|)$, then $(A, |\cdot|)$ has a bounded approximate identity;
\item If $|\cdot|$ is any norm on $A_e$ such that $|\cdot|$ is a Q-norm on $A$ and $A$ is closed in $(A_e, |\cdot|)$, then $|\cdot|$ is a Q-norm on $A_e$;
\item If $|\cdot|$ is a Q-norm on $A_e$, then $A$ is closed in $(A_e, |\cdot|)$.
\end{enumerate}
\end{prop}

\begin{proof}
\emph{(1)} We know that $|\cdot|_{op} \leq |\cdot|$ on $A$. Since $|\cdot|$ is a minimal norm on $A$, it is regular. Hence $|\cdot|_{op} = |\cdot|$.\\
\emph{(2)} Same arguments as statement \emph{(1)} above.\\
\emph{(3)} Let $|\cdot|$ be a minimal norm on $A$. Suppose that there exists a norm $|||\cdot|||$ such that $|||\cdot||| \leq |\cdot|_{op}$. Since $|\cdot|_{op} = |\cdot|$ and $|\cdot|$ is minimal norm on $A$, $|\cdot|_{op} = |\cdot| = |||\cdot|||$ on $A$. Let $a+\lambda1 \in A_e$. Then
\begin{eqnarray*}
|a+\lambda1|_{op} & = & \sup\{|ax+\lambda x|: |x| \leq 1, \; x \in A\}\\
& = &\sup\{|||ax+\lambda x|||: |||x||| \leq 1, \; x \in A\}\\
& = & |||a +\lambda 1|||_{op} \leq |||a + \lambda1|||.
\end{eqnarray*}
Hence $|\cdot|_{op} = |||\cdot|||$ on $A_e$. Therefore $|\cdot|_{op}$ is minimal norm on $A_e$.\\
\emph{(4)} $A$ is dense in $(A_e, |\cdot|)$. Therefore there exists a bounded sequence $e_n$ such that $e_n \longrightarrow 1$. Therefore $e_nx \longrightarrow x$ for every $x \in A$. Hence the proof. \\
\emph{(5)} This proof is similar to the proof of \cite[Theorem 2.7(2)]{DeKa:13}.\\
\emph{(6)} Define $\varphi_{\infty}(a + \lambda 1) = \lambda \; (a+ \lambda 1 \in A_e)$. Since every complex homomorphism is continuous in any Q-norm, $A = \ker{\varphi_{\infty}}$ is closed in $(A_e, |\cdot|)$.
\end{proof}

\begin{defn}
Let $(A, \|\cdot\|)$ be a Banach algebra. The \emph{permanent spectral radius} $r_{pA}(\cdot)$ on $A$ is defined as $r_{pA}(a) = \inf\{|a|: |\cdot| \text{ is a norm on } A\} \, (a \in A)$. Then $A$ has \emph{P-property} if every nonzero closed ideal $I$ in $A$ contains an elements $a$ with $r_{pA}(a) > 0$.
\end{defn}

If $A$ has SEP, then some necessary conditions on $A$ and $A_e$ are listed out in the next result. The terminologies used in the following theorem are standard \cite{Me:92}.

\begin{thm}\label{mayer}
Let $A$ be a non-unital, semisimple Banach algebra. Assume that $A$ has SEP. Then
\begin{enumerate}
\item $A$ has P-property;
\item Every norm on $A$ majorizes a minimal norm;
\item $A_e$ has P-property;
\item Every norm on $A_e$ majorizes a minimal norm;
\item If $|\cdot|$ is a minimal norm on $A_e$, then $(A_e, |\cdot|) \, \widetilde{}$ has P-property.
\end{enumerate}
\end{thm}

\begin{proof}
\emph{(1)} Suppose that $A$ does not have a P-property. Therefore there exists non-zero closed ideal $I$ of $A$ such that $r_{pA}(a) = 0 \; (a \in I)$. Note that $r_{pA}(a) = \inf\{\|a\|: \|\cdot\| \text{ be any norm on A}\}$ (see \cite[Page 79]{Me:92}). Since $A$ has SEP, therefore $r(a) = r_{pA}(a) = 0 \; (a \in I)$. Take non-zero $x \in A$. Then $xa \in I$. Therefore $r(xa) = 0 \; (x \in A)$. Then $a \in radA$. Since $A$ is semisimple, it is contradiction. Hence $A$ has P-property.\\
\emph{(2)} This is \cite[Theorem 2(C)]{Me:92}.\\
\emph{(3)} Let $I$ be a non-zero closed ideal in $A_e$. Therefore $J = I \cap A$ is a non-zero closed ideal in $A$. We know that $r_{pA}(a) \leq |a|$ for any norm $|\cdot|$ on $A_e$. Therefore $r_{pA}(a) \leq r_{pA_e}(a)$. Since $A$ has P-property, there exists $c \in J \setminus \{0\}$ such that $r_{pA}(c) >0$. Therefore $r_{pA_e}(c) \geq r_{pA}(c) >0$. Hence $A_e$ has P-property.\\
\emph{(4)} Since $A_e$ has P-property, the result follows from \cite[Theorem 2(C)]{Me:92}.\\
\emph{(5)} This is \cite[Theorem 2(G)]{Me:92}.
\end{proof}

\begin{exam}
The condition ``semisimplicity'' in Theorem \ref{mayer} can not be omitted. For example, consider $\omega(n)=e^{-n^2}(n \in \mathbb N)$. Then $\ell^1(\mathbb N, \omega)$ is a radical Banach algebra. So clearly $\ell^1(\mathbb N, \omega)$ has SEP but it does not have P-property. For each $k \in \mathbb N$, define $\|f\|_k = \sum_{n=1}^{\infty}|f(n)| e^{-n^2-2nk} \, (f \in \ell^1(\mathbb N, \omega))$. Then each $\|\cdot\|_k$ is a norm on $\ell^1(\mathbb N, \omega)$, $\|\cdot\|_k \geqq \|\cdot\|_l$ but $\|\cdot\|_k \ncong \|\cdot\|_l$ for any $k < l$. Thus the Statements (1) and (2) of Theorem \ref{mayer} do not hold.
\end{exam}

The following is our main result.

\begin{thm}
Let $(A, \|\cdot\|)$ be a non-unital, semisimple Banach algebra. Suppose that $A$ has SEP. Then the following are equivalent.
\begin{enumerate}
\item $A_e$ has SEP;
\item $A$ is closed in $(A_e, |\cdot|)$ for every norm $|\cdot|$ on $A_e$;
\item $0 \in \overline{co}V_A(a)$ with respect to any norm on $A$;
\item $|\cdot|_{op} \cong |\cdot|_1$, for every minimal norm $|\cdot|$ on $A_e$.
\end{enumerate}
\end{thm}

\begin{proof}
\emph{(1)} $\Rightarrow$ \emph{(2)} $A_e$ has SEP. Therefore every norm is Q-norm on $A_e$. Hence the result follows from Proposition \ref{main1}\emph{(6)}.\\
\emph{(2)} $\Rightarrow$ \emph{(3)} Let $|\cdot|$ be any norm on $A$. Since $A$ is semisimple, $|\cdot|_{op}$ is a norm on $A_e$. By hypothesis, $A$ is closed in $(A_e, |\cdot|_{op})$. Therefore the map $\varphi_{\infty}(a+\lambda1) = \lambda$ is $|\cdot|_{op}$-continuous. Note that $||\varphi_{\infty}|| = 1 = \varphi_{\infty}(1)$. Therefore $0 = \varphi_{\infty}(a) \in V_{A_e}^{op}(a;1)$. Hence $0 \in V_{A_e}^{op}(a;1) \; (a \in A)$. By \cite[Theorem 2.4(1)]{DePa:21}, $0 \in \overline{co}V_A(a) \; (a \in A)$.\\
\emph{(3)} $\Rightarrow$ \emph{(4)} Let $|\cdot|$ be any minimal norm on $A_e$. Then by Proposition \ref{main1}\emph{(2)}, $|\cdot|$ is regular on $A$. Hence by Theorem \ref{C}, $|\cdot|_{op} \cong |\cdot|_1$.\\
\emph{(4)} $\Rightarrow$ \emph{(1)} Let $|||\cdot|||$ be any norm on $A_e$. By Proposition \ref{mayer}\emph{(2)}, there exists minimal norm $|\cdot|$ on $A_e$ such that $|\cdot| \leq ||||\cdot|||$ on $A_e$. Therefore by Proposition \ref{main1}\emph{(2)}, $|\cdot| =|\cdot|_{op}$ on $A_e$. By assumption, $|\cdot| = |\cdot|_{op} \cong |\cdot|_1$. Since $A$ has SEP, $|\cdot|_1$ is Q-norm on $A$ and hence on $A_e$. Since $|\cdot| \leq ||||\cdot|||$ on $A_e$. $|||\cdot|||$ is Q-norm on $A_e$. Thus $A_e$ has SEP.
\end{proof}

\begin{conjec}
We strongly believe that if $A$ has SEP, then $A_e$ has SEP at least for semisimple Banach algebras.
\end{conjec}

\noindent
\textbf{Acknowledgement :} {The authors are thankful to the reviewer for reading the manuscript carefully, and giving valuable suggestions and comments for the improvement of the paper. The second author is also  thankful to the Council of Scientific and Industrial Research (CSIR), New Delhi (India), for providing Senior Research Fellowship.}


\begin{thebibliography}{Dillo 83}
\bibitem[BhDe:95]{BhDe:95} S. J. Bhatt, and H. V. Dedania, \emph{Uniqueness of the uniform norm and adjoining identity in Banach algebras}, Proc. Indian Acad. Sci. (Math. Sci.), 105(4)(1995)405-409.
\bibitem[BhDe:96]{BhDe:96} S. J. Bhatt, and H. V. Dedania, \emph{Banach algebras with unique uniform norm}, Proc. American Math. Soc, 124(2)(1996)579-584.
\bibitem[BoDu:71]{BoDu:71} F. F. Bonsall, and J. Duncan, \emph{Numerical Ranges of Operators on Normed Spaces and of Elements of Normed Algebras}, London Math. Soc., Lecture Note Series 2, 1971.
\bibitem[BoDu:73]{BoDu:73} F. F. Bonsall, and J. Duncan, \emph{Complete Normed Algebras}, Springer, Berlin, 1973.
\bibitem[Da:00]{Da:00} H. G. Dales, \emph{Banach Algebras and Atometic Continuity}, Oxford University Press, 2000.
\bibitem[DaDe:09]{DaDe:09} P. A. Dabhi, and H. V. Dedania, \emph{On the uniqueness of uniform norms and C*-norms}, Studia Mathematica, 191(3)(2009)263-270.
\bibitem[DeKa:13]{DeKa:13} H. V. Dedania, and H. J. Kanani, \emph{A non-unital *-algebra has UC*NP if and only if its unitization has UC*NP}, Proc. American Math. Soc, 141(11)(2013)3905-3909.
\bibitem[DPa:21]{DPa:21} H. V. Dedania, and J. G. Patel, \emph{The operator norm on weighted discrete semigroup algebras $\ell^1(S, \omega)$}, Proc. American Math. Soc, 149(12)(2021)5313-5319.
\bibitem[DePa:21]{DePa:21} H. V. Dedania, and A. B. Patel, \emph{Spatial numerical range in non-unital, normed algebras}, Applications and Applied Mathematics: An International Journal, to appear.
\bibitem[GaHu:89]{GaHu:89} A. K. Gaur, and T. Husain, \emph{Spatial numerical ranges of elements of Banach algebras}, Internat. J. Math. \& Math. Sci. 12(4)(1989)633-640.
\bibitem[GaKo:93]{GaKo:93} A. K. Gaur, and Z. V. Kovarik, \emph{Norms on unitizations of Banach algebras}, Proc. American Math. Soc, 117(1)(1993)111-113.
\bibitem[Ka:16]{Ka:16} H. J. Kanani, \emph{Spectral and uniqueness properties in various Banach algebra products}, Ph.D. Thesis, Sardar Patel University, 2016.
\bibitem[Me:91]{Me:91} M. J. Meyer, \emph{The spectral extension property and extension of multiplicative linear functions}, Proc. American Math. Soc, 112(3)(1991)855-861.
\bibitem[Me:92]{Me:92} M. J. Meyer, \emph{Minimal incomplete norms on Banach algebra}, Studia Mathematica, 102(1)(1992)77- 85.
\bibitem[Ta:00]{Ta:00} S. Takahasi, \emph{Spatial numerical ranges of elements of subalgebras of $C_0(X)$}, Internat. J. Math. \& Math. Sci., 23(12)(2000)827-831.
\bibitem[ToYo:89]{ToYo:89} B. J. Tomiuk, and B. Yood, \emph{Incomplete normed algebra norms on Banach algebras}, Studia Mathematica, 95(1989)119-132.
\end{thebibliography}
\end{document}